\documentclass[12pt]{article}
\usepackage{latexsym, amssymb}

\DeclareMathAlphabet\gothic{U}{euf}{m}{n}

\begin{document}

\newcommand{\gotg}{\gothic{g}}

\newcommand{\divv}{\mathop{\rm div}}
\newcommand{\Vol}{\mathop{\rm Vol}}
\newcommand{\Ri}{{\bf R}}
\newcommand{\Ni}{{\bf N}}
\newcommand{\proof}{\mbox{\bf Proof} \hspace{5pt}} 
\newcommand{\ruimte}{\vskip10.0pt plus 4.0pt minus 6.0pt}
\newcommand{\one}{1\hspace{-4.5pt}1}

\newtheorem{lemma}{Lemma}
\newtheorem{thm}[lemma]{Theorem}
\newtheorem{cor}[lemma]{Corollary}
\newtheorem{voorb}[lemma]{Example}
\newtheorem{rem}[lemma]{Remark}
\newtheorem{prop}[lemma]{Proposition}
\newtheorem{stat}[lemma]{{\hspace{-5pt}}}

\newenvironment{remarkn}{\begin{rem} \rm}{\end{rem}}
\newenvironment{exam}{\begin{voorb} \rm}{\end{voorb}}
\newcounter{teller}
\renewcommand{\theteller}{\Roman{teller}}
\newenvironment{tabel}{\begin{list}%
{\rm \bf \Roman{teller}.\hfill}{\usecounter{teller} \leftmargin=1.1cm
\labelwidth=1.1cm \labelsep=0cm \parsep=0cm}
                      }{\end{list}}

\thispagestyle{empty}

\begin{center}
{\Large{\bf {Positivity and strong ellipticity }}}\\[4mm]
{\large A.F.M. ter Elst$^1$, Derek W. Robinson and Yueping Zhu$^2$}\\[4mm]

Centre for Mathematics   and its Applications  \\
Mathematical Sciences Institute \\
 Australian National University  \\
 Canberra, ACT 0200  \\
Australia
\vspace{2cm}

{\large{\bf Abstract}}
\end{center}

\begin{list}{}{\leftmargin=1.8cm \rightmargin=1.8cm}
\item
 We consider partial differential operators $H=-\divv(C\nabla)$
in divergence form on $\Ri^d$ with a positive-semidefinite, symmetric, matrix $C$ 
of real  $L_\infty$-coefficients and  establish that $H$  is strongly elliptic 
if and only if the associated 
 semigroup kernel satisfies local lower bounds, or, if and only if the kernel satisfies
Gaussian upper and lower bounds.
\end{list}

\vfill

\vspace{2cm}
\noindent
September 2004

\vspace{4mm}
\noindent
Keywords: Elliptic operator, semigroup, kernel, upper bounds, lower bounds.

\vspace{1mm}

\noindent
AMS  Classification: 35Jxx.

\vspace{4mm}

\noindent
{\bf Home institutions:}    \\[2mm]
\begin{tabular}{@{}cl@{\hspace{10mm}}cl}
1. & Department of Mathematics  &
  2. & Department of Mathematics \\
& \hspace{15mm} and Computing Science &
  & Nantong  University  \\
& Eindhoven University of Technology &
  & Nantong, 226007   \\
& P.O. Box 513 &
  & Jiangsu Province  \\
& 5600 MB Eindhoven &
  & P.R. China  \\
& The Netherlands &
  & {}
\end{tabular}

\mbox{}

\thispagestyle{empty}

\newpage

\setcounter{page}{1}

The classical Nash--De Giorgi \cite{Nash} \cite{DG} theory analyzes 
positive second-order partial differential operators in divergence form, i.e., operators
\begin{equation}
H=-\sum^d_{i,j=1}\partial_i\,c_{ij}\,\partial_j
\label{eo1.1}
\end{equation}
where $\partial_i=\partial/\partial x_i$, the
 coefficients $c_{ij}$ 
are real $L_\infty$-functions and  
the  matrix $C=(c_{ij})$ is assumed to be symmetric and positive-definite almost-everywhere.
The starting point of the theory is  the strong ellipticity assumption, 
\begin{equation}
C\geq \mu \,I>0
\label{eo1.2}
\end{equation}
almost-everywhere, and 
the principal conclusion is the local H{\"o}lder  continuity of weak solutions
of the associated elliptic and parabolic equations.
In Nash's approach the H{\"o}lder
continuity of the elliptic solution is derived as a 
corollary of continuity of   the parabolic solution and the latter is established 
by an iterative argument from good upper and lower bounds
on the fundamental solution.
Aronson \cite{Aro} subsequently improved Nash's bounds and proved that the 
fundamental solution of the parabolic equation, the heat kernel, satisfies 
Gaussian upper and  lower bounds.
Specifically the kernel $K$  of the semigroup $S$
is a symmetric  function over $\Ri^d\times\Ri^d$ satisfying bounds
\begin{equation}
a'\,G_{b';t}(x-y)\leq K_t(x\,;y)\leq a\,G_{b;t}(x-y)
\label{eo1.3}
\end{equation}
uniformly for all $x,y\in\Ri^d$ and $t>0$ where $G_{b;t}(x)=t^{-d/2}e^{-b|x|^2 t^{-1}}$
and $a,a',b,b'>0$.
Background information and references can be found in the books and reviews \cite{Dav2}
\cite {DER4} \cite{Gri3} \cite{Stroock1} \cite{Stroock2}.

In this note we observe that  a converse statement is true. 
If $H$ is an elliptic operator of the form  (\ref{eo1.1}) then the corresponding heat kernel
satisfies the Aronson bounds (\ref{eo1.3}) if and only if $H$ satisfies the strong ellipticity 
condition (\ref{eo1.2}).
In fact we show that (\ref{eo1.2}) and  (\ref{eo1.3}) are both equivalent to lower bounds
$K_t(x\,;y)\geq a\,t^{-d/2}$ for all $|x-y|\leq rt^{1/2}$ and $t\in\langle0,1]$.

\smallskip

In the Nash--De Giorgi theory
the strong ellipticity assumption (\ref{eo1.2}) is first used to give a precise definition
of  $H$   as a positive self-adjoint operator on $L_2(\Ri^d)$
through quadratic form techniques.
Specifically one  defines the quadratic form $h$ on $L_2(\Ri^d)$ by
\begin{equation}
h(\varphi)=
\sum^d_{i,j=1}\int_{\Ri^d}dx\,
\overline{(\partial_i\varphi)(x)} \, c_{ij}(x) \, (\partial_j\varphi)(x)
\;\;\;.
\label{eo1.4}
\end{equation}
with domain 
$D(h)=W^{1,2}(\Ri^d)=\bigcap^d_{i=1}D(\partial_i)=D(\Delta^{1/2})$
 where  $\Delta$ denotes the self-adjoint
Laplacian, i.e., $\Delta=-\sum^d_{i=1}\partial_i^2$, on $L_2(\Ri^d)$.
Then $h$ is  positive, symmetric, densely-defined and as a direct consequence
of (\ref{eo1.2}) it is also closed.
Therefore there is a unique, positive, self-adjoint operator $H$, with $D(H)\subset D(h)$,
 canonically associated with $h$. 
In particular $(\varphi, H\varphi)=h(\varphi)$ for all $\varphi\in D(H)$.

Since our intention is to analyze the operator $H$ without
the  strong ellipticity assumption (\ref{eo1.2}) the foregoing definition of $H$
is not applicable and one has to adopt an alternative approach.
One can still introduce the form $h$ as above but there is no
reason for the form to be closable.
(For examples of non-closable $h$ see \cite{FOT}, Theorem~3.1.6.)
 We circumvent this problem by a `viscosity' method.

Let  $l$ be the closed quadratic form associated with the Laplacian $\Delta$, i.e.,
\begin{equation}
l(\varphi)=\sum^d_{i=1}\|\partial_i\varphi\|_2^2=\|\Delta^{1/2}\varphi\|_2^2
\label{echse1;1}
\end{equation}
with $D(l)=D(\Delta^{1/2})$.
Then for each  $\varepsilon\in\langle0,1]$ define 
 $h_\varepsilon$ by $D(h_\varepsilon)=D(h)=D(l)$ and
\[
h_\varepsilon(\varphi)=h(\varphi)+\varepsilon\,l(\varphi) 
\]
where $h$ denotes the form given by (\ref{eo1.4}).
Since $h$ is positive the form $h_\varepsilon$ satisfies the strong ellipticity condition
\begin{equation}
h_\varepsilon(\varphi)\geq
\varepsilon\,l(\varphi)
\label{eo2.1}
\end{equation}
for all $\varphi\in D(h)$.
In addition it satisfies the upper bounds
\begin{equation}
h_\varepsilon(\varphi)\leq
(1+\|C\|)\,l(\varphi)
\label{eo2.2}
\end{equation}
 where $\|C\|$ is the essential supremum of the matrix norm of $C(x)=(c_{ij}(x))$.
It follows immediately from (\ref{eo2.1}) and (\ref{eo2.2}) that 
 $h_\varepsilon$ is closed. 
Therefore there is a positive self-adjoint operator
 $H_\varepsilon$ canonically associated with $h_\varepsilon$.
The operator $H_\varepsilon$ is the strongly elliptic operator with coefficients 
$C+\varepsilon I$.
But  $\varepsilon\mapsto h_\varepsilon(\varphi)$ decreases monotonically as $\varepsilon$ decreases
for each $\varphi\in D(h)$.
Therefore it  follows from a result of Kato, \cite{Kat1} Theorem~VIII.3.11, that the $H_\varepsilon$
converge in the strong resolvent sense, as $\varepsilon\to0$, to a positive self-adjoint operator
$H_0$ which we will refer to as the {\bf viscosity operator} with coefficients 
$C=(c_{ij})$.
The strong resolvent convergence also implies that  the positive contractive semigroups 
$S^{(\varepsilon)}$ generated by the 
$H_\varepsilon$ converge strongly to the semigroup $S^{(0)}$ generated by $H_0$.
Therefore $S^{(0)}$ is positive and contractive on $L_2(\Ri^d)$.

Let   $h_0$ denote the form associated with $H_0$,
i.e., $D(h_0)=D(H_0^{1/2})$ and $h_0(\varphi)=\|H_0^{1/2}\varphi\|_2^2$.
There is an alternative method of defining $h_0$ which shows that it has more universal significance.

One may associate with any positive quadratic form $h$ a unique maximal closable minorant $h_r$, i.e.,
$h_r$ is the largest closable positive quadratic form which is majorized by $h$ (see \cite{bSim5} \cite{DalM}).
Then $h_0$ is the closure of $h_r$.
In particular, if $h$ is closable then $h_0$ is its closure.
In addition, $h_0$ is the largest closed positive quadratic form which is majorized by $h$.
Consequently $D(h)\subseteq D(h_0)$.
One may characterize $D(h_0)$ as the vector space of all $\varphi \in L_2$ for which  there are 
$\varphi_1,\varphi_2,\ldots \in D(h)$ such that $\lim_{n \to \infty} \varphi_n = \varphi$ in $L_2$ and 
$\liminf_{n \to \infty} h(\varphi_n) < \infty$.
Moreover, $h_0(\varphi)$ equals the minimum of all $\liminf_{n \to \infty} h(\varphi_n)$, 
where the minimum is taken over all $\varphi_1,\varphi_2,\ldots \in D(h)$ such that 
$\lim_{n \to \infty} \varphi_n = \varphi$ in $L_2$.
(See \cite{bSim4}, Theorem~3.)
In convergence theory $h_0$ is variously called the
lower semi-continuous regularization of $h$ \cite{ET}, page~10, or the  relaxed form 
\cite{DalM}, page~28.

\smallskip

The following theorem  gives a precise formulation of the  characterizations
of strong ellipticity mentioned above.
Other characterizations are given in Proposition~\ref{peo3.1}.

\begin{thm} \label{tchse2}
Let $H_0$ be the viscosity operator with coefficients $C=(c_{ij})$
and $K^{(0)}$ the distribution kernel of the  positive contraction semigroup $S^{(0)}$
generated by $H_0$.
The following conditions are equivalent.
{\samepage
\begin{tabel}
\item\label{tchse2-1}
There is a $\mu>0$ such that   $C\geq\mu \,I$ almost everywhere.
\item\label{tchse2-3}
There are $a,r>0$ such that for all $t\in\langle0,1]$ one has
\[
K^{(0)}_t(x\,;y)\geq a\,t^{-d/2}
\]
for almost every $(x,y) \in \Ri^d \times \Ri^d$ with $|x-y|\leq rt^{1/2}$.
\item\label{tchse2-2}
$K_t^{(0)}$ is a bounded  function satisfying the Aronson Gaussian bounds $(\ref{eo1.3})$.
\end{tabel}}
\end{thm}

The implication \ref{tchse2-1}$\Rightarrow$\ref{tchse2-2} follows by the Nash--Aronson 
estimates and obviously \ref{tchse2-2}$\Rightarrow$\ref{tchse2-3}.
Therefore the proof of the theorem is reduced to establishing that \ref{tchse2-3}$\Rightarrow$\ref{tchse2-1}.
The proof is based on
a variation of an argument of Carlen, Kusuoka and Stroock \cite{CKS}
which requires a  different formulation of  strong ellipticity.

\smallskip

The following proposition gives several related characterizations of
strong ellipticity in terms of the forms $h$, $h_0$ and the corresponding operators.
It is well known (see for example Folland \cite{Fol3}, Theorem~7.17)
that strong ellipticity is equivalent to a G\aa rding inequality
and this may be expressed in terms of either form.

\begin{prop}\label{peo3.1}
Let $H_0$ be the viscosity operator with coefficients $C=(c_{ij})$.
Moreover, let $h$ and~$l$ be the forms given by $(\ref{eo1.4})$ and $(\ref{echse1;1})$ with 
common domain  $\bigcap^d_{i=1}D(\partial_i)$ and let $h_0$ denote the form associated 
with $H_0$.
The following conditions are equivalent.
\begin{tabel}
\item\label{ieo3.0}
The form $h$ is closed.
\item\label{ieo3.0.1}
$h = h_0$.
\item\label{ieo3.1}
There is a $\mu>0$ such that $C\geq\mu \,I$ almost everywhere.
\item\label{ieo3.2.1}
There is a $\mu>0$ such that $h\geq \mu\,l$.
\item\label{ieo3.2.2}
There are $\mu>0$ and $\nu \geq 0$ such that $h\geq \mu\,l - \nu \, I$.
\item\label{ieo3.2}
There is a $\mu>0$ such that $H_0\geq \mu\,\Delta$ in the quadratic form
sense.
\item\label{ieo3.1.2}
There are $\mu>0$ and $\nu\geq 0$ such that $H_0\geq \mu\,\Delta-\nu \,I$ 
in the quadratic form sense.
\end{tabel}

\end{prop}
\proof\
We shall prove that 
\ref{ieo3.1}$\Rightarrow$\ref{ieo3.2.1}$\Rightarrow$\ref{ieo3.0}$\Rightarrow$\ref{ieo3.0.1}$\Rightarrow$\ref{ieo3.1.2}$\Rightarrow$\ref{ieo3.2.2}$\Rightarrow$\ref{ieo3.1}
and 
\ref{ieo3.2.1}$\Rightarrow$\ref{ieo3.2}$\Rightarrow$\ref{ieo3.1.2}.

The implication \ref{ieo3.1}$\Rightarrow$\ref{ieo3.2.1} is trivial.
Since $h \leq \|C\| \, l$ and $l$ is closed the implication \ref{ieo3.2.1}$\Rightarrow$\ref{ieo3.0} is straightforward.
The implication \ref{ieo3.0}$\Rightarrow$\ref{ieo3.0.1} is trivial.

If $h=h_0$ then the vector space $D(\Delta^{1/2}) = D(h) = D(h_0)$ is a Banach space with 
respect to the norm
$\varphi \mapsto (h_0(\varphi) + \|\varphi\|_2^2)^{1/2}$.
But $D(\Delta^{1/2})$ with the graph norm is also a Banach space.
Hence  there is a $c > 0$ such that 
$\|\Delta^{1/2} \varphi\|_2^2 + \|\varphi\|_2^2 \leq c \, (h_0(\varphi) + \|\varphi\|_2^2)$ 
for all $\varphi \in D(\Delta^{1/2})$ as a consequence of the closed graph theorem.
Therefore one deduces \ref{ieo3.1.2}.

The implication \ref{ieo3.1.2}$\Rightarrow$\ref{ieo3.2.2} is evident since $h_0 \leq h$.

The implication \ref{ieo3.2.2}$\Rightarrow$\ref{ieo3.1} follows the proof of Theorem~7.17 in \cite{Fol3}.
Let $\varphi\in C_c^\infty(\Ri^d)$, $k\in\Ri$ and $\xi\in\Ri^d$.
Define $\varphi_k \in D(h) \cap D(l)$ by $\varphi_k(x)=e^{ikx.\xi}\varphi(x)$.
Then one calculates that 
\[
\lim_{k\to\infty}k^{-2} h(\varphi_k) = \int_{\Ri^d}dx\,|\varphi(x)|^2\, (\xi,C(x)\xi)
\;\;\;.  \]
But 
\[
\lim_{k\to\infty}k^{-2} (\mu \, l(\varphi_k) - \nu \,\|\varphi_k\|_2^2)
=\mu \int_{\Ri^d}dx\,|\varphi(x)|^2\,|\xi|^2
\;\;\;.
\]
Since $h(\varphi_k)\geq \mu\, l(\varphi_k) - \nu \|\varphi_k\|_2^2$
one deduces that
\[
\int_{\Ri^d}dx\,|\varphi(x)|^2 \, (\xi,C(x)\xi)
\geq\mu\,\int_{\Ri^d}dx\,|\varphi(x)|^2\,|\xi|^2
\]
Then one concludes that $C\geq\mu I$ almost-everywhere.
This proves the implication \ref{ieo3.2.2}$\Rightarrow$\ref{ieo3.1}.

Next, if \ref{ieo3.2.1} is valid then $\mu \, l$ is a closed positive quadratic form with 
$\mu \, l \leq h$.
Hence $\mu \, l \leq h_0$ and \ref{ieo3.2} is valid.
The implication \ref{ieo3.2}$\Rightarrow$\ref{ieo3.1.2} is obvious.\hfill$\Box$

\ruimte

As a final preliminary to the proof of the missing implication in Theorem~\ref{tchse2}
we need some information on Dirichlet forms \cite{FOT} \cite{BH}.

\smallskip

It is easy to verify that 
$h(|\varphi|) \leq h(\varphi)$ and $h(0 \vee \varphi \wedge \one) \leq h(\varphi)$
for all real valued $\varphi \in D(h)$.
If $\varphi \in D(h_0)$ is real valued then there are $\varphi_1,\varphi_2,\ldots \in D(h)$ such that 
$\lim \varphi_n = \varphi$ in $L_2$ and 
$h_0(\varphi) = \lim h(\varphi_n)$.
Then $\lim |\varphi_n| = |\varphi|$ in $L_2$
and 
$\liminf h(|\varphi_n|) \leq \liminf h(\varphi_n) = h_0(\varphi)$.
So $|\varphi| \in D(h_0)$ and $h_0(|\varphi|) \leq h_0(\varphi)$.
Similarly, $0 \vee \varphi \wedge \one \in D(h_0)$ and 
$h_0(0 \vee \varphi \wedge \one) \leq h_0(\varphi)$.
Therefore $h_0$ is a Dirichlet form  and 
$S^{(0)}$ extends to a positive contraction semigroup 
on all the $L_p$-spaces, which we will also denote by $S^{(0)}$.
It then follows from the positivity and contractivity that the 
semigroup $S^{(0)}$ satisfies
\begin{equation}
0\leq S^{(0)}_t\one\leq \one
\label{eeo2.10}
\end{equation}
for all $t>0$ on $L_\infty(\Ri^d)$.
(In fact one can prove that $S^{(0)}_t\one= \one$
but this is not straightforward (see \cite{ERSZ1}, Proposition~3.6) and it is not necessary in the sequel.)

By the contractivity of $S^{(0)}$ and spectral theory one has 
\[
h_0(\varphi)\geq t^{-1}(\varphi, (I-S^{(0)}_t)\varphi)
\]
for all $\varphi\in D(h_0)$ and $t>0$.
But one deduces from (\ref{eeo2.10}) that 
\[
\|\varphi\|_2^2=(\varphi,\varphi)\geq(S^{(0)}_t\one,|\varphi|^2)=(|\varphi|^2,S^{(0)}_t\one)
\]
for all $t>0$, where $(\,\cdot\,,\,\cdot\,)$ denotes the duality between $L_p$ and $L_q$.
Then it follows from  self-adjointness of $S^{(0)}_t$ and (\ref{eeo2.10}) that
\begin{eqnarray*}
h_0(\varphi)
& \geq & (2t)^{-1}\Big( (S^{(0)}_t\one,|\varphi|^2)+(|\varphi|^2,S^{(0)}_t\one)
    - (\varphi, S^{(0)}_t\varphi) - (S^{(0)}_t\varphi,\varphi) \Big) 
\end{eqnarray*}
for all $\varphi\in D(h_0)$ and $t>0$.
This gives a related estimate in terms of the distribution kernel.

For each $\varphi\in C_c^\infty(\Ri^d)\subset D(h)\subseteq D(h_0)$ choose a $\chi\in C_c^\infty(\Ri^d)$
such that $0\leq \chi\leq1$ and $\chi=1$ on the support of $\varphi$.
Then since $S^{(0)}$ is positive $S^{(0)}_t\chi\leq S^{(0)}_t\one$ and
\begin{eqnarray*}
h_0(\varphi)
& \geq & (2t)^{-1}\Big( (S^{(0)}_t\chi,\chi|\varphi|^2)+(\chi|\varphi|^2,S^{(0)}_t\chi)
    - (\chi\varphi, S^{(0)}_t\chi\varphi) - (S^{(0)}_t\chi\varphi,\chi\varphi) \Big) \\[5pt]
&=& (2t)^{-1} \int_{\Ri^d \times \Ri^d} d(x,y) \,
       K^{(0)}_t(x\,;y)\,\chi(x) \, \chi(y) \, |\varphi(x)-\varphi(y)|^2
\end{eqnarray*}
for all $t>0$.
This is the starting point of the Carlen--Kusuoka--Stroock argument to establish that
\ref{tchse2-3}$\Rightarrow$\ref{tchse2-1} in Theorem~\ref{tchse2}.

\ruimte

\noindent{\bf End of proof of Theorem~\ref{tchse2} \hspace{5pt}}\
Choose a smooth positive function $\rho$ with support in $\langle-r,r\rangle$
such that $\rho\leq 1$ and  $\rho(x) = 1$ for all $x \in \Ri^d$ with if $|x|\leq r/2$.
Then  the previous estimate gives
\[
h_0(\varphi)
\geq (2t)^{-1} \int_{\Ri^d \times \Ri^d} d(x,y) \, 
    K^{(0)}_t(x\,;y)\,\rho(|x-y|^2 t^{-1})\,\chi(x) \, \chi(y) \, 
|\varphi(x)-\varphi(y)|^2
\]
for $\varphi\in C_c^\infty(\Ri^d)$, $t>0$ and $\chi\in C_c^\infty(\Ri^d)$
such that $0\leq \chi\leq1$ and $\chi=1$ on the support of~$\varphi$.
Then it follows from Condition~\ref{tchse2-3} that
\[
h_0(\varphi)
\geq a\,(2t)^{-1}\int_{\Ri^d}dx\int_{\Ri^d}dy\,t^{-d/2}\,\rho(|x-y|^2 t^{-1}) \, 
\chi(x) \, \chi(y) \, |\varphi(x)-\varphi(y)|^2 
\;\;\; .  \]
But the left hand side is independent of the choice of $\chi$ so by the monotone convergence theorem
\[
h_0(\varphi)
\geq a\,(2t)^{-1}\int_{\Ri^d}dx\int_{\Ri^d}dy\,t^{-d/2}\,\rho(|x-y|^2 t^{-1}) \, 
|\varphi(x)-\varphi(y)|^2 
\]
for all $\varphi\in C_c^\infty(\Ri^d)$ and $t\in\langle0,1]$.
Therefore if ${\widehat\varphi}$ denotes the Fourier transform of $\varphi$ then
\begin{eqnarray}
h_0(\varphi)
&\geq& 
a\,t^{-1}\int_{\Ri^d}dx\,t^{-d/2}\,\rho(|x|^2 t^{-1})\int_{\Ri^d}d\xi\, |{\widehat\varphi}(\xi)|^2\,
(1-\cos\xi.x) \nonumber \\[5pt]
&=& 
a\,t^{-1}\int_{\Ri^d}dx\,\rho(|x|^2)\int_{\Ri^d}d\xi\, |{\widehat\varphi}(\xi)|^2 \,
(1-\cos t^{1/2}\xi.x) \nonumber  \\[5pt]
&=&2a\,\int_{\Ri^d}d\xi\, |{\widehat\varphi}(\xi)|^2\,\int_{\Ri^d}dx\,\rho(|x|^2)\,
t^{-1} \sin^2 (2^{-1} t^{1/2} \xi.x) \nonumber 
\end{eqnarray}
for all $\varphi\in C_c^\infty(\Ri^d)$ and $t\in\langle0,1]$.
Thus  in the limit $t\to0$ one has
\[
h_0(\varphi)
\geq 2^{-1}a\, \int_{\Ri^d}d\xi\, |{\widehat\varphi}(\xi)|^2\,\int_{\Ri^d}dx\,\rho(|x|^2)\, (\xi.x)^2
=\mu\int_{\Ri^d} d\xi\, |{\widehat\varphi}(\xi)|^2\, |\xi|^2= \mu\,l(\varphi)
\]
for all $\varphi\in C_c^\infty(\Ri^d)$
with $\mu>0$.
Then since $h_0\leq h$, by the discussion preceding Theorem~\ref{tchse2}, one has
$h(\varphi)\geq \mu\,l(\varphi)$ for all $\varphi\in C_c^\infty(\Ri^d)$.
But  the coefficients of $h$ are bounded and $C_c^\infty(\Ri^d)$ is a core for 
$W^{1,2}(\Ri^d) = D(l) = D(h)$.
So $h(\varphi)\geq \mu\,l(\varphi)$ for all $\varphi \in D(l)$.
Thus Condition~\ref{ieo3.2.1} in Proposition~\ref{peo3.1} is satisfied.
But this is equivalent to Condition~\ref{ieo3.1} of the proposition
which is just a repetition of the strong ellipticity hypothesis, Condition~\ref{tchse2-1} in 
Theorem~\ref{tchse2}.
\hfill$\Box$

\ruimte

Although much work in recent years has been devoted to the derivation of Gaussian
upper bounds on semigroup kernels Theorem~\ref{tchse2} demonstrates that Gaussian 
lower bounds
are in fact the important feature in understanding the general behaviour of the kernels.
The local small time lower bounds in Condition~\ref{tchse2-3} of the theorem encapsulate 
all the information contained in the Aronson upper and lower bounds.
The lower bounds reflect the correct small $t$ behaviour and this is enough to derive
the behaviour of the semigroup and its kernel for all $t$.

It is also interesting to note that in quite general circumstances 
(see, for example, \cite{Cou4}) the Gaussian upper bounds suffice to prove that 
 Gaussian lower bounds are equivalent to H\"older continuity of the kernel.
In particular each of the equivalent conditions of the theorem implies that the semigroup
kernel is H\"older continuous.

\smallskip

It is also possible to extend the theorem to the setting of subelliptic operators on Lie groups.
Let $a_1,\ldots,a_d$ be a vector space basis for the Lie algebra $\gotg$ of a Lie group $G$.
For all $i \in \{ 1,\ldots,d \} $ let $A_i$ be the infinitesimal generator of the 
one parameter group $t \mapsto L(\exp(-t a_i))$, where $L$ is the left regular representation
in $L_2(G)$.
For all $i,j \in \{ 1,\ldots,d \} $ let $c_{ij} \in L_\infty(G)$ and suppose that
the $c_{ij}$ are real and symmetric.
One can define as above a viscosity operator $H_0$
corresponding to the formal expression $- \sum_{i,j=1}^d A_i \, c_{ij} \, A_j$.

Next let $d' \leq d$ and suppose that $a_1,\ldots,a_{d'}$ generate the Lie algebra $\gotg$.
Associated to $a_1,\ldots,a_{d'}$ one can define a modulus $|\cdot|'$ on $G$ and 
a local dimension $D' \in \Ni$, i.e., $\Vol \{ g \in G : |g|' < \rho \} \asymp \rho^{D'}$ for 
$\rho \in \langle0,1]$.
Then one has the following theorem.

\begin{thm} \label{tchse3}
Let $H_0$ be the viscosity operator with coefficients $(c_{ij})$
and $K^{(0)}$ the distribution kernel of the positive contraction semigroup $S^{(0)}$
generated by $H_0$.
The following conditions are equivalent.
\begin{tabel}
\item\label{tchse3-1}
There is a $\mu>0$ such that 
\[
(\varphi, H_0 \varphi) \geq \mu \sum_{i=1}^{d'} \|A_i \varphi\|_2^2
\]
for all $\varphi \in \bigcap_{i=1}^{d'} D(A_i)$.
\item\label{tchse3-3}
There are $a,r>0$ such that for all $t\in\langle0,1]$ one has
\[
K^{(0)}_t(g\,;h) \geq a \, t^{-D'/2}
\]
for almost every $(g,h) \in G \times G$ with $|g h^{-1}|' \leq rt^{1/2}$.

\item\label{tchse3-2}
There are $a,a',b,b',\omega,\omega' > 0$ such that 
\[
a' \, t^{-D'/2} e^{-\omega' t} e^{- b' (|g h^{-1}|')^2 t^{-1}}
\leq K^{(0)}_t(g\,;h)
\leq a \, t^{-D'/2} e^{\omega t} e^{- b (|g h^{-1}|')^2 t^{-1}}
\]
for all $t > 0$ and $g,h \in G$.
\end{tabel}
\end{thm}
The implication \ref{tchse3-1}$\Rightarrow$\ref{tchse3-2} is in \cite{ER18}, 
the implication \ref{tchse3-2}$\Rightarrow$\ref{tchse3-3} is trivial
and the implication \ref{tchse3-3}$\Rightarrow$\ref{tchse3-1} is as in the proof 
of Theorem~\ref{tchse2}, but instead of the scaling of $\rho$ used in the above proof one has to use the 
maps $\gamma_t$ as in \cite{ER13} Section~3.
We omit the technical details.

\ruimte

Finally the situation is quite different for second-order real divergence form operators
which are degenerate \cite{ERSZ1}.
Then the kernel is positive but not necessarily strictly positive even if the operator
is subelliptic. 
One may construct examples for which the kernel vanishes on the loci of degeneracy.
In particular one cannot expect any type of Gaussian lower bound.
Nevertheless subellipticity and a condition of uniform strict positivity suffice 
to deduce Gaussian upper bounds which incorporate the correct large $t$ behaviour
(see \cite{ERSZ1} for details).

\section*{Acknowledgements}

This work was carried out whilst  the first author was visiting the
Centre for Mathematics and its Applications at the ANU.
He wishes to thank the Australian Research Council for support 
and the CMA for  hospitality.
The third author was an  ARC Research Associate for the duration of the collaboration.


\begin{thebibliography}{ERSZ0}

\bibitem[Aro]{Aro}
{\sc Aronson, D.G.}, Bounds for the fundamental solution of a parabolic
  equation.
\newblock {\em Bull.\ Amer. Math.\ Soc.} {\bf 73} (1967),  890--896.

\bibitem[BoH]{BH}
{\sc Bouleau, N. {\rm and} Hirsch, F.}, {\em Dirichlet forms and analysis on
  Wiener space}, vol.\ 14 of de Gruyter Studies in Mathematics.
\newblock Walter de Gruyter \& Co., Berlin, 1991.

\bibitem[CKS]{CKS}
{\sc Carlen, E.A., Kusuoka, S. {\rm and} Stroock, D.W.}, Upper bounds for
  symmetric Markov transition functions.
\newblock {\em Ann.\ Inst.\ Henri Poincar\'e} {\bf 23} (1987),  245--287.

\bibitem[Cou]{Cou4}
{\sc Coulhon, T.}, Off-diagonal heat kernel lower bounds without Poincar\'e.
\newblock {\em J. London Math.\ Soc.} {\bf 68} (2003),  795--816.

\bibitem[{Dal}]{DalM}
{\sc {Dal Maso}, G.}, {\em An introduction to {$\Gamma$}-convergence}, vol.\ 8
  of Progress in Nonlinear Differential Equations and their Applications.
\newblock Birkh{\"a}user Boston Inc., Boston, MA, 1993.

\bibitem[Dav]{Dav2}
{\sc Davies, E.B.}, {\em Heat kernels and spectral theory}.
\newblock Cambridge Tracts in Mathematics 92. Cambridge University Press,
  Cambridge etc., 1989.

\bibitem[DeG]{DG}
{\sc {De Giorgi}, E.}, Sulla differenziabilit\`a e l'analiticit\`a delle
  estremali degli integrali multipli regolari.
\newblock {\em Mem.\ Accad.\ Sci.\ Torino cl.\ Sci.\ Fis.\ Mat.\ Nat.} {\bf 3}
  (1957),  25--43.

\bibitem[DER]{DER4}
{\sc Dungey, N., Elst, A.F.M. ter {\rm and} Robinson, D.W.}, {\em Analysis
  on Lie groups with polynomial growth}, vol.\ 214 of Progress in Mathematics.
\newblock Birkh{\"a}user Boston Inc., Boston, 2003.

\bibitem[EkT]{ET}
{\sc Ekeland, I. {\rm and} Temam, R.}, {\em Convex analysis and variational
  problems}.
\newblock North-Holland Publishing Co., Amsterdam, 1976.

\bibitem[ElR1]{ER13}
{\sc Elst, A.F.M. ter {\rm and} Robinson, D.W.}, Weighted subcoercive
  operators on Lie groups.
\newblock {\em J. Funct.\ Anal.} {\bf 157} (1998),  88--163.

\bibitem[ElR2]{ER18}
\leavevmode\vrule height 2pt depth -1.6pt width 23pt, Second-order subelliptic
  operators on Lie groups II: real measurable principal coefficients.
\newblock In {\sc Balakrishnan, A.~V.}, ed., {\em Proceedings for the First
  International Conference of Semigroups of Operators: Theory and Applications,
  Newport Beach, California}, vol.\ 42 of Progress in nonlinear differential
  equations and their applications. Birkh{\"a}user Verlag, Basel, 2000,
  103--124.

\bibitem[ERSZ]{ERSZ1}
{\sc Elst, A.F.M. ter, Robinson, D.W., Sikora, A. {\rm and} Zhu, Y.},
  Second-order operators with degenerate coefficients, 2004.

\bibitem[Fol]{Fol3}
{\sc Folland, G.B.}, {\em Introduction to partial differential equations}.
\newblock Mathematical Notes 17. Princeton University Press, Princeton, 1976.

\bibitem[FOT]{FOT}
{\sc Fukushima, M., Oshima, Y. {\rm and} Takeda, M.}, {\em Dirichlet forms and
  symmetric Markov processes}, vol.\ 19 of de Gruyter Studies in Mathematics.
\newblock Walter de Gruyter \& Co., Berlin, 1994.

\bibitem[Gri]{Gri3}
{\sc Grigor'yan, A.}, Estimates of heat kernels on Riemannian manifolds.
\newblock In {\em Spectral theory and geometry $($Edinburgh, 1998$)$}, vol.\
  273 of London Math.\ Soc.\ Lecture Note Ser.,  140--225. Cambridge Univ.\
  Press, Cambridge, 1999.

\bibitem[Kat]{Kat1}
{\sc Kato, T.}, {\em Perturbation theory for linear operators}.
\newblock Second edition, Grundlehren der mathematischen Wissenschaften 132.
  Springer-Verlag, Berlin etc., 1984.

\bibitem[Nas]{Nash}
{\sc Nash, J.}, Continuity of solutions of parabolic and elliptic equations.
\newblock {\em Amer.\ J. Math.} {\bf 80} (1958),  931--954.

\bibitem[Sim1]{bSim4}
{\sc Simon, B.}, 
Lower semicontinuity of positive quadratic forms.
\newblock {\em Proc.\ Roy.\ Soc.\ Edinburgh Sect.\ A} {\bf 79} (1977),
  267--273.

\bibitem[Sim2]{bSim5}
\leavevmode\vrule height 2pt depth -1.6pt width 23pt, 
A canonical decomposition for quadratic forms with
  applications to monotone convergence theorems.
\newblock {\em J.\ Funct.\ Anal.} {\bf 28} (1978),  377--385.

\bibitem[Str1]{Stroock1}
{\sc Stroock, D.W.}, Diffusion semigroups corresponding to uniformly elliptic
  divergence form operators.
\newblock In {\sc Az{\'e}ma, J., Meyer, P.A. {\rm and} Yor, M.}, eds., {\em
  S{\'e}minaire de probabilit{\'e}s XXII}, Lecture Notes in Mathematics 1321.
  Springer-Verlag, Berlin etc., 1988,  316--347.

\bibitem[Str2]{Stroock2}
\leavevmode\vrule height 2pt depth -1.6pt width 23pt, 
Estimates for the heat kernel of second order elliptic
  operators.
\newblock In {\em Nonlinear partial differential equations and their
  applications. Coll\`ege de France Seminar, Vol.\ XII (Paris, 1991--1993)},
  vol.\ 302 of Pitman Res.\ Notes Math.\ Ser.,  226--235. Longman Sci.\ Tech.,
  Harlow, 1994.

\end{thebibliography}
\end{document}